\documentclass[9pt]{amsart}

\usepackage[latin1]{inputenc}
\usepackage[T1]{fontenc}
\usepackage{enumerate}
\usepackage{diagbox}  

\usepackage[english]{babel}

\usepackage{amsmath,amsfonts,amssymb,amsthm,amscd}

\newtheorem{theorem}{Theorem}[section]
\newtheorem{lemma}[theorem]{Lemma}
\newtheorem{corollary}[theorem]{Corollary}

\newtheorem{claim}{Claim}

\theoremstyle{definition}
\newtheorem{definition}[theorem]{Definition}

\newtheorem{problem}{Problem}
\newtheorem{question}{Question}

\theoremstyle{remark}
\newtheorem{remark}[theorem]{Remark}

\numberwithin{equation}{section}

\newcommand{\N}{\mathbb N}

\title[Supercyclic composition operators on the space of smooth functions]{Supercyclic weighted composition operators on the space of smooth functions}
\author[J.\ B\`{e}s and C.\ Foster]{J.\ B\`{e}s and C.\ Foster}
\address{J.~B\`{e}s, Department of Mathematics and Statistics,
Bowling Green State University,
Bowling Green, Ohio 43403,
USA.}
\email{jbes@bgsu.edu}
\address{C. Foster, Department of Mathematics and Statistics,
Bowling Green State University,
Bowling Green, Ohio 43403,
USA.}
\email{fostercd@bgsu.edu}
\date{July 22, 2025}

\subjclass{Primary 47A16; Secondary 47B33, 47B38}

\keywords{Weighted composition operators, chaotic operators, hypercyclic operators, supercyclic operators, mixing operators}

\allowdisplaybreaks[3]
\begin{document}
\begin{abstract}
A weighted composition operator on the space of scalar-valued smooth functions on an open set $\Omega$ of d-dimensional Euclidean space is supercyclic if and only if it is weakly mixing, and it is strongly supercyclic if and only if it is mixing. Every mixing such operator is chaotic.  In the one-dimensional case, it is supercyclic if and only if it is mixing and if and only if it is chaotic.
\end{abstract}
\maketitle
%{\large   
\section{Introduction}
%upper case only
\setcounter{equation}{0} 
\numberwithin{equation}{section} 
The study of dynamical properties of weighted composition operators on function spaces reached several advances over the past decade (see e.g. \cite{AlbaneseJordaMele2022,BeltranJordaMurillo2019,golprz2021,CA2017,Kalmes2017,Przestacki2017,Zajac,ZhangZhou2016}). We are interested here in the work by Przestacki \cite{Przestacki2017} and by Albanese, Jord\'a and Mele \cite{AlbaneseJordaMele2022}
on the dynamics of this class of operators on the space of smooth functions.
Throughout this paper, $\mathbb{K}$ denotes the real or complex scalar field and $C^{\infty} (\Omega, \mathbb{K})$ the space of $\mathbb{K}$-valued smooth functions on an open subset $\Omega$ of $\mathbb{R}^d$. 
Przestacki~\cite{Przestacki2017} completely determined characterizations of when a weighted composition operator
\[
C_{\omega, \psi}:C^{\infty} (\Omega, \mathbb{K})\to C^{\infty} (\Omega, \mathbb{K}), f(\cdot)\mapsto \omega(\cdot) (f\circ\psi)(\cdot)
\] is hypercyclic, weakly mixing, and mixing (see Subsection~{1.1} for definitions and notation). Here $\omega :\Omega\to \mathbb{K}$ and $\psi:\Omega\to \Omega$ are smooth. 
To state his characterizations precisely, we recall that the symbol $\psi$ is said to be {\bf run-away} (respectively, {\bf strongly run-away}) provided for each compact subset $K$ of $\Omega$ we have
\[\psi_n (K) \cap K = \emptyset\]
for some $n\in\mathbb{N}$ (respectively, for all large $n\in \mathbb{N}$). Here $\psi_n$ denotes the $n$-fold composition of $\psi$ with itself.

\begin{theorem} \label{T1.1} {\rm (Przestacki \cite[Theorem~{3.5} and Theorem~{3.6}]{Przestacki2017})}

 Let $\Omega \subset \mathbb{R}^{d}$ be open, $d\in \mathbb{N}$, and let $\psi : \Omega \to \Omega$ and $\omega : \Omega \to \mathbb{K}$ be smooth. The following are equivalent:
\begin{enumerate}
\item[(1)]\ The operator $C_{\omega, \psi}$ is hypercyclic on $C^{\infty} (\Omega, \mathbb{K})$.

\item[(2)]\ The operator $C_{\omega, \psi}$ is weak mixing on $C^{\infty} (\Omega, \mathbb{K})$.

\item[(3)]\ The following conditions are satisfied: 
\begin{enumerate}
\item[(a)]\ For every $x \in \Omega$, we have $\omega (x) \neq 0$.

\item[(b)] $\psi$ is injective.

\item[(c)]\ For every $x \in \Omega$, we have $det[\psi' (x)] \neq 0$.

\item[(d)]\ $\psi$ has the run-away property.
\end{enumerate}
\end{enumerate}
Moreover, $C_{\omega, \psi}$ is mixing on $C^{\infty} (\Omega, \mathbb{K})$ if and only if it satisfies conditions (a), (b) and (c) above and $\psi$ has the strong run-away property.
\end{theorem}

Przestacki also showed that when $d=1$ and $\Omega  =\mathbb{R}$ the properties of hypercyclicity, mixing and chaos are equivalent within this class of operators, and recently Albanese et al~\cite{AlbaneseJordaMele2022} showed that in this case supercyclicity is also an equivalent property.
\begin{theorem} \label{T:1.3}  {\rm (Przestacki (\cite[Theorem 4.2]{Przestacki2017}), Albanese et al~(\cite[Theorem~3.5(iii)]{AlbaneseJordaMele2022}))                           } 
 $C_{\omega, \psi}$ is supercyclic on $C^{\infty}(\mathbb{R}, \mathbb{K})$ if and only if it is mixing, and if and only if it is chaotic.
\end{theorem}

In summary, we have with the above results full characterizations for hypercyclicity, weak-mixing and mixing of weighted composition operators on $C^{\infty} (\Omega, \mathbb{K})$, and when $\Omega=\mathbb{R}$ we also have characterizations for supercyclicity and chaos. Motivated by this we consider the following.

\begin{problem} \label{Q1}
Is there a characterization for supercyclicity of weighted composition operators on  $C^{\infty} (\Omega, \mathbb{K})$ for other cases  than when $\Omega = \mathbb{R}$? 
\end{problem}
\begin{problem}\label{Q2}
Is there a characterization for chaos of weighted composition operators on  $C^{\infty} (\Omega, \mathbb{K})$ for other cases than when $\Omega = \mathbb{R}$? 
\end{problem}

We show in this paper that for an arbitrary open set $\Omega\subseteq \mathbb{R}^d$ a weighted composition operator is supercyclic on $C^{\infty} (\Omega, \mathbb{K})$ if and only if it is weakly mixing, and that it is strongly supercyclic precisely when it is mixing (Theorem~\ref{T:2.1}). In particular, supercyclicity has the same characterization provided by Prestacki's Theorem~\ref{T1.1}.
We also show that in this general setting every mixing weighted composition operator is chaotic (Theorem~\ref{T:4.2}).
Finally, we show that in the one dimensional case  the equivalences of Theorem~\ref{T:1.3} also hold when $\Omega$ is a proper open subset of $\mathbb{R}$  (Theorem~\ref{T:2.2}).

We  conclude the introduction with notation and preliminaries. We show    Theorem~\ref{T:2.1}  in Section 2,   Theorem~\ref{T:4.2} 
in Section 3, and Theorem~\ref{T:2.2} in Section 4.

\subsection{Notation and Preliminaries}

Recall that the space $C^{\infty} (\Omega, \mathbb{K})$ is a separable Fr\'echet space when equipped with the  the family of seminorms 
\[
\{ \| . \|_{K, n}: \ K\subset \Omega \mbox{ compact}, n\in \mathbb{N}_0 \},\]
where 
\[ \lvert \lvert f \rvert \rvert_{K, n} =\max\limits_{\lvert \alpha \rvert \leq n} \max\limits_{x \in K}  \lvert 
 \frac{\partial^{\lvert \alpha 
\rvert} f}{\partial x^{\alpha}} (x) \rvert \ \ \ \ \ (f\in C^{\infty} (\Omega, \mathbb{K})).\] 
Here $\mathbb{N}_0$ denotes the set of non-negative integers,  and 
$
\frac{\partial^{\lvert \alpha \rvert} }{\partial x^{\alpha}} 
$
the operator of partial differentiation of order $|\alpha |=   \alpha_1+\dots +\alpha_d$ that is induced by   the multi-index $\alpha=(\alpha_1,\dots \alpha_d)\in \mathbb{N}_0^d$. We also let $C^p(\Omega, \mathbb{K})$ $(p\in\mathbb{N})$ denote the Fr\'echet space of scalar valued functions whose partial derivatives of order not larger than $p$ are all continuous, which is topologized by the family of seminorms
\[
\{ \| . \|_{K, n}: \ K\subset \Omega \mbox{ compact, $n\in \mathbb{N}_0$ with $n\le p$}                           \} .\] Finally, we let $C^p(\Omega, \Omega)$ denote the set of self-maps $F=(f_1,\dots, f_d)$ of $\Omega$ for which $f_j\in C^p(\Omega, \mathbb{R})$ for each $j=1,\dots, d$. For $F\in C^p(\Omega, \Omega)$ and $x\in\Omega$ we denote the derivative of $F$ at $x$ by $F'(x)$ which is a linear map from $\mathbb{R}^d$ to $\mathbb{R}^d$ and which we may identify with its Jacobian matrix denoted by $[F'(x)]$. Similarly, the derivative of any $f\in C^p(\Omega, \mathbb{C})$ at $x$, denoted by $f'(x)$, is a linear map from $\mathbb{R}^d$ to $\mathbb{R}^2$.

\begin{definition} 
An operator $T$ on a separable Fr\'echet space $X$ is said to be {\bf hypercyclic }(respectively, {\bf supercyclic}) provided there exists a vector $g$ in $X$ whose orbit
\[
\mbox{Orb}(g, T)=\{ g, Tg, T^2g,\dots \}
\]
(respectively, whose projective orbit
$
\mathbb{K} \mbox{Orb}(g, T) = \{ \lambda T^ng: \ n\in\mathbb{N}_0, \lambda\in\mathbb{K} \} 
$) 
is dense in $X$. Such  $g$ is called a {\bf hypercyclic vector} (respectively, a {\bf supercyclic vector}) for $T$.  Also, we say that $T$ is {\bf hereditarily hypercyclic } (respectively, {\bf hereditarily supercyclic}) with respect to a given strictly increasing sequence $(n_k)$ of positive integers provided for each subsequence $(n_{k_j})$ of $(n_k)$ there exists some $g$ in $X$ for which
\[
\{ g, T^{n_{k_1}}g, T^{n_{k_2}}g, \dots \}      \ \ \ \mbox{ (respectively, for which } \mathbb{K} \{ g, T^{n_{k_1}}g, T^{n_{k_2}}g, \dots \} \mbox{ ) }
\]
is dense in $X$. Finally, $T$ is {\bf strongly hypercyclic} (respectively, {\bf strongly supercyclic}) provided it is hereditarily hypercyclic (respectively, hereditarily supercyclic) with respect to the full sequence $(n)$ of positive integers.

\end{definition}

When $X$ is a separable Fr\'echet space the properties of being hypercyclic, hereditarily hypercyclic, and strongly hypercyclic are equivalent to being (topologically) transitive, weak-mixing and mixing, respectively, while chaos in Devaney's sense is equivalent to being both transitive and having a dense set of periodic points \cite{godefroy_shapiro1991operators}.
\begin{definition}
An operator $T$ on a topological vector space $X$ is said to be {\bf transitive} (respectively, {\bf mixing}) provided for each non-empty open subsets $U$, $V$ of $X$ we have
\[
T^n(U) \cap V \ne \emptyset
\]
for some $n$ (respectively, for all large $n$). Also, it is said to be {\bf weakly mixing} provided its direct sum $(x,y)\overset{T\oplus T}{\mapsto} (Tx,Ty)$ is a transitive operator on $X\times X$. Finally $T$ is said to be {\bf chaotic} provided it is transitive and its set
$
\cup_{n\in\mathbb{N}} \mbox{Ker}(T^n-I)
$
of periodic points is dense in $X$.
\end{definition}

We note that a supercyclic operator on an infinite-dimensional Fr\'echet space must have dense range. Indeed, this holds under the weaker assumption of being supercyclic with respect to the weak topology.

\begin{lemma} \label{L:denserange}
Let $X$ be a separable infinite dimensional Fr\'echet space over a real or complex scalar field $\mathbb{K}$. If $T:X\to X$ is weakly supercyclic, then it has dense range.
\end{lemma}
\begin{proof}
Let $h\in X$ be a weakly supercyclic vector for $X$. Since $\mathbb{K} \cdot \mbox{Orb}(h, T)\subseteq \mbox{span}\{ h \} \cup T(X)$, we have
\[
X=\overline{\mathbb{K} \cdot \mbox{Orb}(h, T)}^w \subseteq \overline{\mbox{span}\{ h \}}^w \cup \overline{T(X)}^w=\mbox{span}\{ h \} \cup \overline{T(X)},
\]
where for a subset $A$ of $X$ we denote here by $\overline{A}$ and $\overline{A}^\omega$ the closure and the weak closure of $A$ in $X$ respectively.
Thus $X\setminus \mbox{span}\{ h \} \subseteq \overline{T(X)}$. Now Baire's Category theorem ensures that one of the two subspaces $\mbox{span}\{ h \}$ or $T(X)$ to be somewhere dense in $X$. But proper closed subspaces of a Fr\'echet space have empty interior, so the assumption that $X$ is infinite-dimensional forces the subspace $T(X)$ to have dense range.
\end{proof}

We refer to \cite{bayart_matheron2009dynamics} and \cite{grosse-erdmann_peris-manguillot2011linear} for general background on linear dynamics.

\section{Supercyclicity equals weak mixing, strong supercyclicity equals mixing}
We show in this section the following.

\begin{theorem} \label{T:2.1} Let $\Omega \subset \mathbb{R}^{d}$ be open, and let $\psi \in C^{p} (\Omega, \Omega)$ and $\omega \in C^{p} (\Omega, \mathbb{K})$, where $p\in\mathbb{N}\cup\{\infty\}$ and $d\in\mathbb{N}$. Then the following are equivalent:
\begin{enumerate}
\item[(1)]\ $C_{\omega, \psi}$ is supercyclic on $C^{p} (\Omega, \mathbb{K})$.
\item[(2)]\ $C_{\omega, \psi}$ is weakly mixing on $C^{p} (\Omega, \mathbb{K})$.
\item[(3)]\ The following conditions are satisfied: 
\begin{enumerate}
\item[(i)]\ For every $x \in \Omega$, we have $\omega (x) \neq 0$.
\item[(ii)]\ $\psi$ is injective.
\item[(iii)]\ For every $x \in \Omega$, we have det$[\psi' (x)] \neq 0$.
\item[(iv)]\ $\psi$ is run-away on $\Omega$.
\end{enumerate}
\end{enumerate}
\noindent Moreover, $C_{\omega, \psi}$ is strongly supercyclic if and only if $C_{\omega, \psi}$ is mixing, and if and only if both $\psi$ is strongly run-away on $\Omega$ and conditions (3) (i) - (iii) hold.
\end{theorem}

We first derive with  Lemma \ref{L:supercyclic,necessary} below some necessary conditions for a weighted composition operator to be supercyclic on $C^{\infty}(\Omega, \mathbb{K})$ and use these to derive with Lemma \ref{L:supercyclic,run-away} the pending necessary condition by which altogether suffice to establish Theorem \ref{T:2.1}.

\begin{lemma} \label{L:supercyclic,necessary}

Let $\Omega \subset \mathbb{R}^{d}$ be open, and let $C_{\omega, \psi}$ be supercyclic on $C^{p} (\Omega, \mathbb{K})$, where $p\in \mathbb{N}\cup\{ 0, \infty \}$. Then the following conditions are satisfied: 
\begin{enumerate}
\item[(a)]\ The multiplier $\omega$ is zero-free. 
\item[(b)]\ The symbol $\psi$ is injective.
\item[(c)]\ If $p>0$, then det$[\psi' (x)] \neq 0$ for each $x\in\Omega$.
\item[(d)]\ If $p>0$, the symbol $\psi$ has no periodic points. 
\end{enumerate}
\end{lemma}

\begin{proof}  We recall that $C_{\omega, \psi}$ must have dense range by Lemma~\ref{L:denserange}.

(a)  Suppose there exists $x_0 \in \Omega$ with $\omega (x_0) = 0$. Then for each $g \in C^p (\Omega, \mathbb{K})$ we have 
\[ C_{\omega, \psi} (g) (x_0) =  \omega(x_0)  (g \circ \psi) (x_0) = 0,\]
forcing the range of $C_{\omega, \psi}$ to be contained in the closed hyperplane
\[ \{f \in C^{p} (\Omega, \mathbb{K}): f(x_0) = 0\}\]
of $C^{p} (\Omega, \mathbb{K})$, contradicting that $C_{\omega, \psi}$ has dense range. 

(b) Suppose that $\psi (x_1) = \psi (x_2)$ for some $x_1, x_2 \in \Omega$ with $x_1 \neq x_2$.  By (a) we know that $\omega (x) \neq 0$ for every $x \in \Omega$. So for each $g \in C^{p} (\Omega, \mathbb{K})$ we have
\[ C_{\omega, \psi} (g) (x_1) = \omega (x_1)  g(\psi (x_1)) = \frac{\omega (x_1)}{\omega (x_2)}  \omega(x_2)  g(\psi (x_2)) = \frac{\omega (x_1)}{\omega (x_2)} C_{\omega, \psi} (g) (x_2),\]
forcing the range of $C_{\omega,\psi}$ to be contained in the closed hyperplane
\[ \{f \in C^{p} (\Omega, \mathbb{K}): f(x_1) = \frac{\omega (x_1)}{\omega (x_2)} \cdot f(x_2)\}\]
of $C^{p} (\Omega, \mathbb{K})$ which contradicts that  $C_{\omega, \psi}$ has dense range. So $\psi$ must be injective.

(c) Suppose there exists $x_0 \in \Omega$ with $\mbox{det}[\psi' (x_0)] = 0$. Then there exists a non-zero vector $h \in \mathbb{R}^{d}$ such that 
\[\psi' (x_0) h = (0, ..., 0),\]
and since $p\ge 1$ for each $g \in C^{p} (\Omega, \mathbb{K})$ we have
\[
\begin{aligned}
C_{\omega, \psi} (g)' (x_0) h & =  (\omega \cdot  (g \circ \psi))' (x_0) h &\\
& = g(\psi (x_0)) \omega'(x_0) h +  \omega (x_0)  g'(\psi (x_0)) \psi' (x_0) h &\\
& =  g(\psi (x_0)) \omega' (x_0) h &\\ 
& =  \omega (x_0)  g(\psi (x_0)) \frac{\omega' (x_0)}{\omega (x_0)} h &\\
& =   C_{\omega, \psi} (g) (x_0)  \frac{\omega' (x_0)}{\omega (x_0)} h.
\end{aligned}
\]
Hence it follows that the range of  $C_{\omega,\psi}$ is contained in the proper closed linear subspace
 \[ \left\{f \in C^{p} (\Omega, \mathbb{K}): f' (x_0) h = f(x_0)  \frac{\omega' (x_0)}{\omega (x_0)} h \right \}\]
  of $C^{p} (\Omega, \mathbb{K})$, which contradicts that
   $C_ {\omega, \psi}$ has dense range. Thus for each $x \in \Omega$ the derivative of $\psi$ at x is nonsingular.

(d) Notice that by Ansari's theorem~\cite[Theorem~2]{Ansari} (which holds on arbitrary topological vector spaces) %%% The SC version of Ansari may be deducted from \cite[Exercise 3.6, p73]{bayart_matheron2009dynamics} %%%
$C_{\omega, \psi}$ is supercyclic if and only if any iterate $C_{\omega, \psi}^r=C_{\tilde{\omega}, \psi_r}$ ($r\in\mathbb{N}$) is supercyclic, where $\tilde{\omega}=\prod_{j=0}^{r-1} \omega\circ \psi_j$. So it suffices to show that $C_{\omega, \psi}$ has no fixed point. Now, by means of contradiction suppose there exists $x \in \Omega$ with $\psi(x) = x$. Throughout this proof we view $\psi'(x)$ as an operator acting on $\mathbb{C}^{d}$ regardless if $\mathbb{K}$ is the real or complex scalar field. Let $\lambda_0$ be an eigenvalue of $\psi'(x)$, and let  $h_0 \in \mathbb{C}^{d}$ be an eigenvector of $\psi'(x)$ with eigenvalue $\lambda_0$. 
Since $C_{\omega, \psi}$ is supercyclic, the commuting family $\mathcal{F} = \{\lambda C^{n}_{\omega, \psi}:n \in \mathbb{N}\cup\{0\}, \lambda \in \mathbb{K} \}$ is topologically transitive on the Baire space $C^{p} (\Omega, \mathbb{K})$ and  the set $SC(C_{\omega, \psi})$ of supercyclic vectors for $C_{\omega, \psi}$,  which coincides with the set of universal vectors for $\mathcal{F}$, is a co-meager set in $C^{p}(\Omega, \mathbb{K})$, see \cite[p27]{bayart_matheron2009dynamics}. So there exists $f \in SC(C_{\omega, \psi})$ satisfying $f(x) \neq 0$ and $f'(x) h_0 \neq 0$. 

Pick $g \in C^{p}(\Omega, \mathbb{K})$ so that $g(x) \neq 0$, $g'(x) h_0 \neq 0$, and $g \notin \bigcup_{n=0}^{\infty} \text{span}\{C^{(n)}_{\omega, \psi} (f)\}$. Replacing g by $g + \epsilon$ for some small $\epsilon \in \mathbb{R}$ if necessary, without loss of generality we may further assume that 
\[
\frac{\lvert g'(x) h_0 \rvert}{\lvert g(x) \rvert} \neq \frac{\lvert f' (x) h_0 \rvert}{\lvert f(x) \rvert}.
\]

\noindent Since $f \in SC(C_{\omega, \psi})$, there exist sequences $(\lambda_l)_{l \in \mathbb{N}}$ and $(n_l)_{l \in \mathbb{N}}$  in $\mathbb{K}$ and $\mathbb{N}$ respectively so that 
\begin{equation} \label{eq:232}
 \lambda_l C^{n_l}_{\omega, \psi} (f)  \underset{\ell\to\infty}{\to} g \hspace{.2cm} \text{in} \hspace{.2cm} C^{p} (\Omega, \mathbb{K}).
\end{equation}
\noindent In particular, since $p\ge 1$ all first order partial derivatives of the above sequence converge uniformly on compact sets to the corresponding partial derivative of the limit. Hence
\begin{equation} \label{eq:233}
(\lambda_l C^{n_l}_{\omega, \psi} (f))' (x)  \underset{\ell\to\infty}{\to} g' (x) \hspace{.2cm} \text{in} \hspace{.2cm} L(\mathbb{C}^{d}, \mathbb{C}^{d}).
\end{equation}

\noindent \textbf{Claim}: For each $n \in \mathbb{N}$ and $F \in C^{p} (\Omega, \mathbb{K})$ we have 
\begin{equation}\label{eq:claim}
(C^{n}_{\omega, \psi} (F))' (x) = F(x) \omega (x)^{n-1} \omega' (x) + \omega (x) (C^{n-1}_{\omega, \psi} (F))'(x) \psi'(x).
\end{equation}
Indeed, by the Product Rule and the Chain Rule
\[
\begin{aligned}
(C_{\omega, \psi}(F))'(x) & = F(\psi(x)) \omega'(x) + \omega(x) F'(\psi(x)) \psi'(x) &\\ 
& = F(x) \omega'(x) + \omega(x) F'(x) \psi'(x) &\\
& = F(x) \omega (x)^{n-1} \omega' (x) + \omega (x) (C^{n-1}_{\omega, \psi} (F))'(x) \psi'(x),
\end{aligned}
\]
so $\eqref{eq:claim}$ holds for $n=1$. Inductively, if $\eqref{eq:claim}$  holds for the case $n=k$ we have
\[
\begin{aligned}
(C^{k+1}_{\omega, \psi}(F))'(x) & = (C^{k}_{\omega, \psi} (C_{\omega, \psi} (F)))' (x) &\\ 
& = (C_{\omega, \psi} (F)) (x) \omega (x)^{k-1} \omega'(x) + \omega(x) (C^{k-1}_{\omega, \psi} (C_{\omega, \psi} (F)))' (x) \psi'(x) &\\ 
& = F(x) \omega(x)^{k}  \omega '(x) + \omega(x) (C^{k}_{\omega, \psi} (F))' (x) \psi' (x),  
\end{aligned}
\]
so it holds for the case $n=k+1$ as well.  Hence the claim holds. Now, letting $F = f$ in the claim, we get 
\begin{equation} \label{eq:234}
\begin{aligned}&(C^{n}_{\omega, \psi} (f))' (x)  = \omega(x)^{n-1} f(x) \omega '(x) + \omega(x) (C^{n-1}_{\omega, \psi} (f))' (x) \psi' (x) \\
& = \omega(x)^{n-1} f(x) \omega '(x) + \omega(x)[\omega (x)^{n-2} f(x) \omega' (x) + \omega(x)(C^{n-2}_{\omega, \psi} (f))'(x) \psi'(x)]\psi'(x) \\
& = 2 \omega(x)^{n-1} f(x) \omega'(x) + \omega(x)^{2}(C^{n-2}_{\omega, \psi} (f))' (x) [\psi' (x)]^{2} \\
& \ \,  \vdots \\
& = n \omega(x)^{n-1} f(x) \omega'(x) + \omega(x)^{n} f' (x) [\psi' (x)]^{n} \\
& = \frac{n}{\omega(x)} C^{n}_{\omega, \psi} (f) (x) \omega' (x) + \frac{1}{f(x)} C^{n}_{\omega, \psi} (f) (x) f'(x) [\psi' (x)]^{n} \\
& = C^{n}_{\omega, \psi} (f) (x) \left( \frac{n}{\omega(x)} \omega'(x) + \frac{1}{f(x)} f'(x) [\psi'(x)]^{n}\right).
\end{aligned}
\end{equation}

\noindent So it follows from $\eqref{eq:234}$ that for all large $l$ we have
\begin{equation} \label{eq:235}
\begin{aligned} \frac{n_l}{\omega(x)} \omega'(x) h_0 &+ \frac{\lambda_0 ^{n_l}}{f(x)} f'(x) h_0  = \left (\frac{n_l}{\omega(x)} \omega'(x) + \frac{1}{f(x)} f'(x) [\psi'(x)]^{n_l}\right) h_0 \\
& = \frac{1}{\lambda_l C^{n_l}_{\omega, \psi} (f)(x)} \lambda_l C^{n_l}_{\omega, \psi} (f)(x) \left(\frac{n_l}{\omega(x)} \omega'(x) + \frac{1}{f(x)} f'(x) [\psi'(x)]^{n_l}\right)h_0 \\
& = \frac{1}{\lambda_l C^{n_l}_{\omega, \psi} (f)(x)} (\lambda_l C^{n_l}_{\omega, \psi} (f))' (x) h_0,
\end{aligned}
\end{equation}
and thus by  $\eqref{eq:232}$ and $\eqref{eq:233}$ we have
\begin{equation} \label{eq:236}
     \frac{n_l}{\omega(x)} \omega'(x) h_0 + \frac{\lambda_0 ^{n_l}}{f(x)} f'(x) h_0 \underset{\ell\to\infty}{\to} \frac{1}{g(x)} g'(x) h_0.
\end{equation}
We now derive a contradiction using $\eqref{eq:236}$. 
We have two cases.

\noindent Case 1: $\lvert \lambda_0 \rvert > 1$. We have 
\[
\begin{aligned}
 \left| n_l \frac{1}{\omega(x)}  \omega'(x) h_0 + \frac{\lambda^{n_l}_0}{f(x)} f'(x) h_0 \right|  & \geq \left| \frac{\lambda^{n_l}_0}{f(x)} f'(x) h_0 \right| - n_l \left|  \frac{1}{\omega(x)}  \omega'(x) h_0 \right| &\\
& = \left| \lambda_0 \right|^{n_l} \left(\frac{\left| f' (x) h_0 \right|}{\left| f(x) \right| } - \frac{n_l}{\left| \lambda_0 \right| ^{n_l}}
\left| \frac{1}{\omega(x)} \omega'(x) h_0 \right| \right) \underset{\ell\to\infty}{\to} \infty,    
\end{aligned}
\]
a contradiction with $\eqref{eq:236}$.

\noindent Case 2: $\lvert \lambda_0 \rvert \leq 1$. If $\omega'(x) h_0 \neq 0$, we have 
\begin{flalign*} \left| n_l \frac{1}{\omega(x)} \omega'(x) h_0 + \frac{\lambda^{n_l}_0}{f(x)} f'(x) h_0 \right| & \geq \left| n_l \frac{1}{\omega(x)} \omega'(x) h_0 \right| - \left| \frac{\lambda^{n_l}_0}{f(x)} f'(x) h_0 \right| &\\
& \geq n_l \frac{\lvert \omega'(x) h_0 \rvert}{\lvert \omega (x) \rvert} - \frac{\lvert f'(x) h_0 \rvert}{\lvert f(x) \rvert}  \underset{\ell\to\infty}{\to} \infty,
\end{flalign*}
again contradicting $\eqref{eq:236}$. Now assume $\omega'(x) h_0 = 0$. Then we have
\begin{equation} \left| n_l \frac{1}{\omega(x)} \omega'(x) h_0 + \frac{\lambda^{n_l}_0}{f(x)} f'(x) h_0 \right| = \frac{\lvert \lambda_0 \rvert ^{n_l}}{\lvert f(x) \rvert} \lvert f'(x) h_0 \rvert \to \begin{cases}
   \ \ \  0 & \mbox{if} \hspace{2mm} \lvert \lambda_0 \rvert < 1 \\
    \frac{\lvert f'(x) h_0 \rvert}{\lvert f(x)\rvert} & \mbox{if} \hspace{2mm} \lvert \lambda_0 \rvert = 1.
\end{cases}
\end{equation} 
But $\frac{\lvert g'(x) h_0 \rvert}{\lvert g(x) \rvert} \neq \frac{\lvert f' (x) h_0 \rvert}{\lvert f(x) \rvert}$ and $\frac{g'(x)h_0}{g(x)}\ne 0$ by our selection of $g$. So in either case we get a contradiction with $\eqref{eq:236}$, and the proof of $(4)$ is now complete. 

\end{proof}

\begin{remark} \label{R:3.2}
By Lemma~\ref{L:denserange}, the same argument shows that conclusions (a)-(c) of Lemma~\ref{L:supercyclic,necessary}  hold under the weaker assumption that $C_{\omega, \psi}$ be weakly supercyclic. Also, conclusions (a) and (b) hold when $p=0$, too. We don't know whether conclusion (d) holds when $p=0$ or under the weaker assumption that $C_{\omega, \psi}$ is weakly supercyclic.
\end{remark}

We say that $\psi:\Omega\to\Omega$ is run-away with respect to a given strictly increasing sequence $(n_j)$ in $\mathbb{N}$ provided for each compact subset $K$ of $\Omega$ there exists $j\in\mathbb{N}$ so that 
\[\psi_{n_j}(K)\cap K=\emptyset.\]
Also, recall that a family $\mathcal{F}=\{ T_\alpha \}_{\alpha\in \Lambda}$ of operators on a topological vector space $X$ is said to be {\em universal} provided there exists a vector $g$ in $X$ for which
\[
\mathcal{F}g=\{ T_\alpha g: \alpha\in \Lambda \}
\]
is dense in $X$. Any such  $g$ is called a {\em universal vector } for $\mathcal{F}$, and when each $T_\alpha$ has dense range and $T_\alpha T_\beta=T_\beta T_\alpha$ for each $\alpha,\beta\in\Lambda$ the family $\mathcal{F}$ is universal if and only if it has a dense set of universal vectors.

\begin{lemma} \label{L:supercyclic,run-away}
Let $p\in \mathbb{N}\cup \{ \infty \}$, let $\omega \in  C^{p} (\Omega, \mathbb{K})$, and let $\psi\in C^{p} (\Omega, \Omega)$. If $C_{\omega, \psi}$ is supercyclic on $C^{p} (\Omega, \mathbb{K})$ with respect to some strictly increasing sequence $(n_j)$ in $\mathbb{N}$, the following hold:
\begin{enumerate}
\item[(i)]\ For each $x\in \Omega$, the set $\overline{\{\psi_{n_j} (x): j \geq 0 \}}^\Omega$ is not compact in $\Omega$.

\item[(ii)]\ The symbol $\psi$ is run-away with respect to $(n_j)$.
\end{enumerate}
\end{lemma}

\begin{proof}
(i) Suppose that for some $x_0 \in \Omega$ the set 
\[
K := \overline{\{\psi_{n_j }(x_0) : j \geq 0\}}^\Omega
\]
is a compact subset of $\Omega$. By Lemma~\ref{L:supercyclic,necessary}(a) and the continuity of $\omega$  there exist positive scalars $M, m > 0$ such that
\[
m < \lvert \omega (x) \rvert < M
\]
 for every $x \in K\cup \psi(K)$.  Since the set of universal vectors for \[
 \mathcal{F}:=\{ \lambda C^m_{\omega, \psi}:  \lambda \in \mathbb{K},   m\in\{ n_j: \ j\in \mathbb{N} \} \cup \{ 0 \}  \} 
 \]
 is dense in $C^{p} (\Omega, \mathbb{K})$, we may pick a universal vector h for $\mathcal{F}$ which satisfies
\begin{equation}\label{eq:1<h<2}
1 < \lvert h(x) \rvert < 2
\end{equation}
 for every $x \in K\cup \psi(K)$. In particular, for each $j\in\mathbb{N}$ and $0\ne \lambda\in\mathbb{K}$ we must have 
\begin{equation}\label{eq:-1}
\left| \frac{\lambda C^{n_j}_{\omega, \psi}(h) (\psi(x_0))}{\lambda C^{n_j}_{\omega, \psi}(h) (x_0)} \right| = \left| \frac{\omega (\psi_{n_j} (x_0)) h(\psi_{n_j+1} (x_0))}{\omega (x_0) h(\psi_{n_j} (x_0))} \right| < \frac{2M}{m}.
\end{equation}
Now, since $C^{p} (\Omega, \mathbb{K})$ is dense in $C(\Omega, \mathbb{K})$, any universal vector for  $\mathcal{F}$ acting on $C^{p} (\Omega, \mathbb{K})$ is also a universal vector for  $\mathcal{F}$ acting on $C(\Omega, \mathbb{K})$. So without loss of generality for the remainder of the proof we may consider $C_{\omega, \psi}$ and $\mathcal{F}$ as acting on $C(\Omega, \mathbb{K})$. 

Notice that $\psi(x_0)\ne x_0$ by Lemma~\ref{L:supercyclic,necessary}(d), so the open subset
\[ U := \left \{f \in C (\Omega, \mathbb{K}): f(x) \neq 0 \, \text{for every} \, x \in K \text{ and} \, \left \lvert \frac{f(\psi (x_0))}{f(x_0)} \right \rvert > \frac{2M}{m}\right \}
\]
of  $C(\Omega, \mathbb{K})$ is non-empty. Also, by $\eqref{eq:-1}$ we have $\lambda C^{n_j}_{\omega, \psi}(h) \not\in U$ for each $j\in\mathbb{N}$ and $\lambda\ne 0$. Hence the universality of $h$  forces that $\lambda_0 h\in U$ for some $0\ne \lambda_0\in \mathbb{K}$ and thus the contradiction
\[
\frac{2M}{m}<\left| \frac{\lambda h(\psi(x_0))}{\lambda h(x_0)}\right| = \left| \frac{ h(\psi(x_0))}{ h(x_0)}\right| <2,
\]
where the last inequality holds by $\eqref{eq:1<h<2}$.

(ii) Assume to the contrary that $\psi$ is not run-away with respect to $(n_j)$. Then there exists a compact set $K \subset \Omega$ such that 
\[
\psi_{n_j} (K) \cap K \neq \emptyset
\] for every $j \in \mathbb{N}$. For each  $j \in \mathbb{N}$ let $x_j \in K$ such that $\psi_{n_j} (x_j) \in K$. 
Arguing as in part (i) there exist  $M, m > 0$ and $h\in C^{p} (\Omega, \mathbb{K})$ universal for $\mathcal{F}$ so that
\begin{equation} \label{eq:FI}
 m < \lvert \omega (x) \rvert < M \ \ \mbox{ and } \ \ 1< |h(x)| <2 
 \end{equation}
for every $x \in K\cup \psi(K)$, and we may also assume $\mathcal{F}$ to be acting on $C(\Omega, \mathbb{K})$.
Again,  for each $j\in\mathbb{N}$ and non-zero scalar $\lambda$ we have 
\begin{equation}\label{eq:0}
\left| \frac{\lambda C^{n_j}_{\omega, \psi}(h) (\psi(x_j))}{\lambda C^{n_j}_{\omega, \psi}(h) (x_j)}\right|= \left| \frac{\omega (\psi_{n_j} (x_j)) h(\psi_{n_j+1} (x_j))}{\omega (x_j) h(\psi_{n_j}(x_j))} \right|< \frac{2M}{m}.
\end{equation}
It now suffices to verify that the open subset
\[
U: = \left\{f \in C(\Omega, \mathbb{K}): f(x) \neq 0 \, \text{and} \, \left| \frac{f(\psi(x))}{f(x)} \right| > \frac{2M}{m} \, \text{for every} \, x \in K \right\}
\]
of  $C(\Omega, \mathbb{K})$ is non-empty. Indeed, since   $\eqref{eq:0}$ forces $\lambda C^{n_j}_{\omega, \psi}(h) \not\in U$ for each $j \in\mathbb{N}$ and $\lambda\ne 0$ we must now have $\lambda_0 h\in U$ for some $\lambda_0\ne 0$, which by $\eqref{eq:FI}$  gives the contradiction
\[
\frac{2M}{m}<\left| \frac{\lambda_0 h(\psi(x_j))}{\lambda_0 h(x_j)}\right| = \left| \frac{ h(\psi(x_j))}{ h(x_j)}\right| <2.
\]
Finally, that $U$ is non-empty was established in \cite[p.1103]{Przestacki2017}: Pick any $d \in C^{\infty} (\Omega, \mathbb{K})$ satisfying $d=0$ on $K$ and  $0<d\le 1$ on $\Omega\setminus K$, and consider
\[
D:=\sum_{n=0}^\infty \frac{1}{(C+1)^{n}} d\circ \psi_n,
\]
where  $C = \frac{2M}{m}$. The uniform convergence of this series ensures that $D\in C (\Omega, \mathbb{K})$. Also Part (i) ensures that for any given $x\in\Omega$ the set $\overline{\{\psi_{n_j} (x): j \geq 0 \}}^\Omega$ is not compact in $\Omega$, forcing that $\psi_{n_j}(x)\notin K$
for some $j\in \mathbb{N}$ and thus that $D(x)>0$. Finally, for each $x\in K$ we have
\[
D(\psi(x)) = \sum_{i=0}^{\infty} \frac{1}{(C+1)^{i}} d(\psi_{i+1} (x)) = (C+1)D(x) - (C+1)d(x) = (C+1)D(x),
\]
so
\[
\frac{D(\psi (x))}{D(x)} = C+1 
\]
and thus  $D \in U$. 
\end{proof}

\begin{remark}
If conclusion (d) of Lemma~\ref{L:supercyclic,necessary} holds when $p=0$, then we may include $p=0$ in the assumptions of Lemma~\ref{L:supercyclic,run-away}. Indeed, when $C_{\omega, \psi}$ is supercyclic with respect to $(n_j)$, Condition (d) of Lemma~\ref{L:supercyclic,necessary} is equivalent to condition (i) of Lemma~\ref{L:supercyclic,run-away}, and to condition (ii) of Lemma~\ref{L:supercyclic,run-away}, regardless of which $p\in \mathbb{N}\cup \{ 0, \infty \}$ we consider.
\end{remark}

\begin{corollary}\label{C:stronglyscimpliesstronglyrunaway} Let $p\in\mathbb{N}\cup\{\infty\}$.
If $C_{\omega,\psi}$ is hereditarily supercyclic on $C^p(\Omega, \mathbb{K})$ with respect to the full sequence $(n)$, then $\psi$ is strongly run-away.
\end{corollary}

We are now ready to show Theorem \ref{T:2.1}.%

\begin{proof}[Proof of Theorem \ref{T:2.1}] The implication $(1) \Rightarrow (3)$  follows by Lemma~\ref{L:supercyclic,necessary} and Lemma~\ref{L:supercyclic,run-away}, and $(2) \Rightarrow (1)$ is immediate.  $(3) \Rightarrow (2)$ Notice that $C_{\omega, \psi}:C^{p}(\Omega, \mathbb{K})\to C^{p}(\Omega, \mathbb{K})$ is quasi-conjugate to $C_{\omega, \psi}:C^{\infty}(\Omega, \mathbb{K})\to C^{\infty}(\Omega, \mathbb{K})$, and the latter is weak-mixing by Theorem~\ref{T1.1}.
For the equivalence in the last statement of the conclusion, notice that any mixing operator on a separable Fr\'echet space is hereditarily supercyclic with respect to the full sequence $(n)$. Conversely, assume $C_{\omega, \psi}$ is hereditarily supercyclic with respect to $(n)$.  To see that it is mixing, by Theorem~\ref{T1.1} it suffices to verify that conditions (3) (i) - (iii) hold and that  $\psi$ is strongly run-away. Conditions (3) (i) - (iii) now follow by the first part of Theorem \ref{T:2.1}, already established. Finally, Corollary~\ref{C:stronglyscimpliesstronglyrunaway} ensures that  $\psi$ is strongly run-away.
\end{proof}

\section{Mixing implies chaos}
\begin{theorem} \label{T:4.2}
Let $\Omega \subset \mathbb{R}^{d}$ be open, and let $\psi \in C^{p} (\Omega, \Omega)$ and $\omega \in C^{p} (\Omega, \mathbb{K})$, where $p\in\mathbb{N}\cup\{\infty\}$. If $C_{\omega, \psi }$ is mixing on $C^p(\Omega, \mathbb{K})$ then it is chaotic.
\end{theorem}

\begin{definition} 
Given $\Omega\subseteq \mathbb{R}^d$  open and non-empty and $\psi:\Omega\to \Omega$, we define
$\psi_0:\Omega\to \Omega$ as the identity self-map of $\Omega$ and for each $n\in\mathbb{N}$ we let $
\psi_n:\Omega\to \Omega$ denote the $n$-fold composition of $\psi$ with itself. If $\psi$ is an injective open mapping, for each $n\in\mathbb{N}$ $\psi_n(\Omega)$ is open and non-empty and we let
$\psi_{-n}:\psi_n(\Omega)\to \Omega$ be defined by the rule
\[
y=\psi_n(x) \mbox{ if and only if } x=\psi_{-n}(y)
\]
\end{definition}
\begin{remark} \label{R:1.9}
Notice that if $p\in\mathbb{N}\cup\{ \infty \}$ and if $\psi\in C^p(\Omega, \Omega)$ is a 1-1 open mapping then for each $n\in \mathbb{N}$ we have 
\[
\psi_n\in C^p(\Omega, \Omega) \ \mbox{ and } \ 
\psi_{-n}\in C^p(\psi_n(\Omega), \Omega).
\]
Also, for each $m, n\in \mathbb{N}\cup\{ 0 \}$ we have
\begin{equation} \label{eq:R1.1}
\begin{aligned}
\psi_n\circ \psi_m&=\psi_{n+m}\ \ =\ \psi_m\circ \psi_n \ \ \ \ \ \ \ \ \ \  \ \ \mbox{ on $\Omega$ }& \ \\
\psi_{-n}\circ \psi_{-m} &=\psi_{-m-n}= \psi_{-m}\circ \psi_{-n}  \ \ \ \ \ \  \mbox{ on $\psi_{m+n}(\Omega)$. }& \ 
\end{aligned}
\end{equation}

Moreover, if in addition  $n\le m$, then we have
\begin{equation}\label{eq:R1.2}
\begin{aligned}
\psi_{-n}\circ \psi_m&=\psi_{m-n}\ \ \ \ \ \ \ \ \ \ \ \  \ \ \ \ \mbox{ on $\Omega$ }\\
\psi_{m}\circ \psi_{-n} &=\psi_{m-n}\ \ \ \ \ \ \ \ \ \ \  \ \ \ \ \mbox{ on $\psi_{n}(\Omega)$ }\\
\psi_{-m}\circ \psi_n&=\psi_{-m+n}\ \ \ \ \ \ \ \ \ \ \ \ \mbox{ on $\psi_{m-n}(\Omega)$ }\\
\psi_{n}\circ \psi_{-m} &=\psi_{-m+n}\ \ \ \ \ \  \ \ \ \ \ \ \mbox{ on $\psi_{m}(\Omega)$. }
\end{aligned}
\end{equation}
\end{remark}

\begin{proof}
It suffices to show $C_{\omega, \psi}$ has a dense set of periodic points. So let $f\in C^p(\Omega, \mathbb{K})$, and let $m\le p$, $\epsilon >0$ and $K\subset \Omega$ compact be given. We want to find $g\in C^p(\Omega, \mathbb{K})$  and $N\in\mathbb{N}$ so that
\begin{equation}.  \label{eq:a1}
 \begin{aligned} 
\| g - f \|_{K, m}&<\epsilon \\
C_{\omega, \psi}^Ng&=g.
\end{aligned}
\end{equation}
Let $0<\rho <\mbox{dist}\{ K, \partial\Omega \}$, and consider the compact subset
\[
L:=\{ x\in \Omega: \ \mbox{dist}(x, K)\le \rho \}
\]
of $\Omega$. The selection of $\rho$ ensures that for each $x\in\Omega$ we have
\[
\begin{aligned}
 x\in \mbox{int}(L) &\Leftrightarrow \ \mbox{dist}(x, K)<\rho \\
 x\in \partial L &\Leftrightarrow \ \mbox{dist}(x, K)=\rho, 
\end{aligned}
\]
and
\[
\partial L=\partial \mbox{int}(L)=L\setminus \mbox{int}(L).
\]
Since $\psi$ is strongly run-away, there exists $N\in \mathbb{N}$ so that for each $n\ge N$ we have
\begin{equation} \label{eq:eqUno}
\psi_n(L)\cap L =\emptyset.
\end{equation}

For each integer $k$ we define
\[
L_k:=\begin{cases}
\ \ \psi_{kN}(L) &\mbox{ if } k\ge 0 \\
\ \ \psi_{-kN}^{-1} ( L) &\mbox{ if } k<0.
\end{cases}
\]
So $L_k$ is always non-empty when $k\ge 0$, while when $k<0$ we have
\[
L_k=\begin{cases} \ \ \ \ \ \  \emptyset
 &\mbox{ if } \ \psi_{-kN}(\Omega)\cap L =\emptyset \\
\psi_{kN} (\psi_{-kN}(\Omega)\cap L) &\mbox{ if } \ \psi_{-kN}(\Omega)\cap L \ne \emptyset. 
\end{cases}
\]
We will use  the following observations.

\begin{remark}  \label{R:FACT1}
\

\begin{enumerate}
\item[(i)]\ 
For each $k\in \mathbb{Z}$ and $y\in L_k$ there exists a unique $x\in L$ so that $y=\psi_{kN}(x)$, and when $k<0$ we also have that $x\in \psi_{-kN}(\Omega)\cap L$.  
\item[(ii)]\ Each non-empty $L_k$ is diffeomorphic to $L$. Indeed,
when $L_k\ne\emptyset$ we have
\[
\psi_N(L_k)=L_{k+1},
\]
and $\psi_N$ maps a neighborhood of $L_k$ onto some neighbourhood of $L_{k+1}$, so for each neighbourhhod $V_{k+1}$ of $\partial L_{k+1}$ there exists a neighbourhood $V_k$ of $\partial L_k$ so that $\psi_N(V_k)\subset V_{k+1}$.
\item[(iii)]\  For any $y\in \Omega$ we have $\psi_N(y)\in \cup_{k\in\mathbb{Z}} L_k$  if and only if  $y\in  \cup_{k\in\mathbb{Z}} L_k$.
\end{enumerate}
\end{remark}

We also use the following.

\begin{claim}
\

\begin{enumerate}
\item[(i)]\ For each $k\in \mathbb{Z}$ with $L_k\ne \emptyset$ we have 
$
\mbox{dist}(L_k, \cup_{j\in\mathbb{Z}\setminus\{ k\} } L_j)>0.
$
\item[(ii)]\  $\mbox{int} (\cup_{k\in\mathbb{Z}} L_k) = \cup_{k\in\mathbb{Z}} \mbox{int}(L_k)$

\item[(iii)]\ The set $\cup_{k\in\mathbb{Z}} L_k$ is closed, and 
$
\partial  (\cup_{k\in \mathbb{Z}}L_k)  =\cup_{k\in\mathbb{Z}} \partial L_k.
$
\end{enumerate}
\end{claim}

Assume Claim~1 holds. Pick some $\varphi\in C^\infty (\Omega, \mathbb{R})$ with $\varphi (\Omega)=[0,1]$ so that
\[
\begin{cases} 
\varphi =1 \ &\mbox{ on } \{ x\in \Omega: \ \mbox{dist}(x, K) \le \frac{1}{3}\rho \} \\
\varphi =0 &\mbox{ on } \Omega\setminus U,
\end{cases}
\]
where $U:= \{ x\in \Omega: \ \mbox{dist}(x, K) < \frac{2}{3}\rho \}$. Then $\tilde{f}:=\varphi f\in C^p(\Omega, \mathbb{K})$ satisfies that
\[
\| \tilde{f}-f \|_{K, m}=0 \ \mbox{ and } \ \ \ \mbox{supp}(\tilde{f})\subset U.
\]
Let $H:=\prod_{j=0}^{N-1} \omega\circ\psi_j$, and
consider the function $g:\Omega\to \mathbb{K}$ given by
\begin{equation} \label{eq:g}
g(y):=\begin{cases} 
\ \left( \frac{ \tilde{f}\circ \psi_{-nN}}{\prod_{j=-n}^{-1} (H\circ \psi_{jN})}\right ) (y) \ &\mbox{ if  $y\in L_n$ for some $n\in\N$,}\\
\ \ \ \ \tilde{f}(y) \ &\mbox{ if $y\in L_0=L$,}\\
\ \ \tilde{f}(\psi_{nN}(y)) \prod_{s=0}^{n-1} H(\psi_{sN}  (y)) \ &\mbox{ if $y\in L_{-n}$ for some $n\in \N$,} \\
\ \ \ \  \ 0 \ &\mbox{ if  $y\in \Omega\setminus \cup_{k\in \mathbb{Z}}L_k$,.} 
\end{cases}
\end{equation}

Notice that $g=\tilde{f}$ on $\mbox{int}(L_0)$ so $g_{\restriction_{\mbox{int}(L_0)}}\in C^p(\mbox{int}(L_0), \mathbb{K})$. Also, since $H\in C^p(\Omega, \mathbb{K})$ and is zero-free, it follows by Remark~\ref{R:FACT1} that $g_{\restriction_{\mbox{int}(L_k)}}\in C^p(\mbox{int}(L_k), \mathbb{K})$ whenever $L_k\ne\emptyset$, and thus by Claim~1(ii) that 
\[
g_{\restriction_V} \in C^p(V , \mathbb{K})
\]
where $V=\mbox{int}(\cup_{k\in\mathbb{Z}} L_k)$. Also $g=0$ on $\Omega\setminus \cup_{k\in \mathbb{Z}} L_k$, so
\[
g_{\restriction_{\Omega\setminus \cup_{k\in \mathbb{Z}}} L_k}\in C^p(\Omega\setminus \cup_{k\in \mathbb{Z}} L_k, \mathbb{K}).
\]
Finally, notice that by $\eqref{eq:g}$ and Claim~1(i) we must have $g=0$ on some neighborhood of $L_0$, since 
 $\tilde{f}=0$ on a neighborhood of $\partial L_0$. Indeed, by the latter and Remark~\ref{R:FACT1}(ii) for each $k\in\mathbb{Z}\setminus\{ 0\}$ with $L_k\ne\emptyset$ we have $\tilde{f}\circ \psi_{-kN} =0$ on some neighborhood of $\partial L_k$, which forces by Claim~1(iii) that $g=0$ on some neighborhood of $\partial \cup_{k\in \mathbb{Z}} L_k$. So $g\in C^p(\Omega, \mathbb{K})$, and it suffices now to show that \[
C_{\omega, \psi}^Ng=g.\]
 Notice first by $\eqref{eq:g}$ and Remark~\ref{R:FACT1} for each $y\in \Omega\setminus  \cup_{k\in\mathbb{Z}} L_k$
\[
g(y)=0=g(\psi_N(y))=H(y) g(\psi_N(y))=C_{\omega, \psi}^N(g)(y).
\]
Also, for $y\in L_0$ we have
\[
C_{\omega, \psi}^Ng(y)=H(y) g(\psi_N(y)) = H(y)  \left( \frac{\tilde{f}\circ \psi_{-N}}{H\circ \psi_{-N}} \right) (\psi_N(y)) =\tilde{f}(y)=g(y).
\]
For $n\in\mathbb{N}$ and $y\in L_n$ we have $\psi_N(y)\in L_{n+1}$ and by $\eqref{eq:g}$
\[
\begin{aligned}
C_{\omega, \psi}^Ng(y)=&H(y) g(\psi_N(y)) \\
&= H(y)  
\left( 
\frac{ \tilde{f}\circ \psi_{-(n+1)N} }{\prod_{j=-(n+1)}^{-1} H\circ \psi_{jN} } 
\right) 
(\psi_N(y))\\
&= \left( \frac{ \tilde{f}\circ \psi_{-nN}}{\prod_{j=-n}^{-1} H\circ \psi_{jN}} \right)(y) =g(y).
\end{aligned}
\]
Finally, if $y \in \cup_{n\in\mathbb{N}} L_{-n}$, we have two cases:
\noindent
Case 1: $y\in L_{-1}$.  Since $\psi_N(y)\in L_0$ by $\eqref{eq:g}$ we have
\[
g(y) =\tilde{f}(\psi_N(y)) H(y) = g(\psi_N(y)) H(y) = C_{\omega. \psi}^Ng(y).
\]
\noindent
Case 2: $y\in L_{-n}$ for some $n\ge 2$.  In this case $\psi_N(y)\in L_{-n+1}$ and $-n+1\le -1$, so
\[
\begin{aligned}
g(y)&= \tilde{f}(\psi_{nN}(y)) \prod_{s=0}^{n-1} H(\psi_{sN}(y)) \\
&=\left( \tilde{f}\circ \psi_{(n-1)N}\right) (\psi_N(y)) H(y) \prod_{s=0}^{n-2} (H\circ \psi_{sN})(\psi_N(y))\\
&=H(y) g(\psi_N(y))\\
&=C_{\omega, \psi}^Ng (y).
\end{aligned}
\]
So $C_{\omega, \psi}^Ng=g$, and it suffices to show Claim~1.

\begin{proof}[Proof of Claim~1]
(i) We first show that 
\begin{equation}\label{eq:ij0}
L_i\cap L_j=\emptyset 
\end{equation}
whenever $i\ne j$. Indeed, interchanging $i$ and $j$ if necessary we may assume that $j=i+n$ for some $n\in\mathbb{N}$.
If $y\in L_i\cap L_j$, there exist $x_i, x_j\in L$ so that 
\begin{equation} \label{eq:ij}
y=\psi_{iN}(x_i)=\psi_{jN}(x_j)
\end{equation}
If $i\ge 0$ we have by $\eqref{eq:R1.1}$ that $\psi_{(i+n)N}=\psi_{iN}\circ \psi_{nN}$ on $\Omega$. So
\[
\psi_{jN}(x_j)=\psi_{(i+n)N}(x_j)=\psi_{iN} (\psi_{nN}(x_j))
\]
 and by $\eqref{eq:ij}$ and the injectivity of $\psi$ we must have
\[
x_i=\psi_{nN}(x_j)\in L\cap \psi_{nN}(L),
\]
a contradiction with $\eqref{eq:eqUno}$. So $i<0$. Now, notice that by $\eqref{eq:R1.2}$ we have
\[
\psi_{iN}\circ \psi_{nN}=\psi_{(i+n)N} \mbox{ on }
\begin{cases} \Omega & \mbox{ if } j=n+i\ge 0\\
\psi_{-j}(\Omega) &\mbox{ if } j=n+i<0.
\end{cases}
\]
Notice also that if $j<0$ then $x_j\in L\cap \psi_{-j}(\Omega)$, so in either case we have
\[
\psi_{jN}(x_j)=\psi_{(i+n)N} (x_j) = \psi_{iN}(\psi_{nN}(x_j))
\]
 and again by $\eqref{eq:ij}$ and the injectivity of $\psi$ we have
\[
x_i=\psi_{nN}(x_j)\in L\cap \psi_{nN}(L),
\]
a contradiction with $\eqref{eq:eqUno}$. So $\eqref{eq:ij0}$ holds. Now,  let $k\in\mathbb{Z}$ be fixed. 
Consider the compact set
\[
\widehat{L_k}:=\{ x\in \Omega: \ \mbox{dist}(x, L_k) \le \frac{1}{3}\rho \}.
\]
Since $\psi$ is strongly run-away on $\Omega$, there exists $N_1$ so that
\begin{equation}\label{eq:hat}
\psi_n(\widehat{L_k}) \cap \widehat{L_k} = \emptyset
\end{equation}
for each $n\ge N_1$. In particular, for each $n\ge N_1$ we have
\[
\mbox{dist}(L_k, L_{k+n})=\mbox{dist} (L_k, \psi_{nN}(L_k))\ge \mbox{dist}(L_k, \psi_{nN}(\widehat{L_k}))\ge  \frac{1}{3}\rho.
\]
We have two cases. Case 1: There exists $j_0$ (necessarily negative) so that $L_j=\emptyset$ for each $j\le j_0$ and $L_{j}\ne \emptyset$ for each $j_0<j$.
By $\eqref{eq:ij0}$ the compact sets $L_k$ and 
$
\cup_{j=j_0,  j\ne k}^{N_1+k} L_j
$
are disjoint and hence at a positive distance, so
\[
\mbox{dist}(L_k, \cup_{j\in\mathbb{Z}\setminus\{ k\} } L_j) \ge \mbox{min}\{ \mbox{dist}(L_k, \cup_{j=j_0,  j\ne k}^{N_1+k} L_j), \frac{1}{3}\rho \} >0.
\]
\noindent
Case 2: $L_j\ne \emptyset$ for each $j\in \mathbb{Z}$. By $\eqref{eq:hat}$ for each $n\ge N_1$ we have
\[
 \widehat{L_k}\cap \psi_{nN}^{-1} (\widehat{L_k}) = \emptyset, 
\]
which by $\eqref{eq:R1.2}$ gives
\[
\begin{aligned}
\frac{1}{3}\rho &\le \mbox{dist}(L_k, \psi_{nN}^{-1} (\widehat{L_k}) )\\
&\le  \mbox{dist}(L_k, \psi_{nN}^{-1} (L_k) )\\
&=\mbox{dist}(L_k, L_{k-n}).
\end{aligned}
\]
and again by $\eqref{eq:ij0}$ we have $L_k$ at a positive distance from $\cup_{j=-N_1+k,  j\ne k}^{N_1+k} L_j$ so
\[
\mbox{dist}(L_k, \cup_{j\in\mathbb{Z}\setminus\{ k\} } L_j) \ge \mbox{min}\{ \mbox{dist}(L_k, \cup_{j=-N_1+k,  j\ne k}^{N_1+k} L_j), \frac{1}{3}\rho \} >0.
\] 

(ii) Clearly the inclusion $\mbox{int} (\cup_{k\in\mathbb{Z}} L_k) \supseteq \cup_{k\in\mathbb{Z}} \mbox{int}(L_k)$ holds. Now, given $x\in \mbox{int} (\cup_{k\in\mathbb{Z}} L_k) $ let $\epsilon >0$ so that 
\[
B(x,\epsilon)\subset \cup_{j\in\mathbb{Z}} L_j
\]
and pick $k\in\mathbb{Z}$ so that $x\in L_k$. By $(i)$ there exists $\delta_k>0$ so that
\[
\mbox{dist}(L_k, \cup_{j\in\mathbb{Z}\setminus\{ k\} } L_j)\ge \delta_k.
\]
Then \[
B(x, \mbox{min}\{ \epsilon, \delta_k \})\subset \cup_{j\in\mathbb{Z}}L_j \setminus \cup_{j\in\mathbb{Z}\setminus\{ k\} } L_j=L_k.\]
So $x\in \mbox{int}(L_k)$, and $(ii)$ holds.

(iii) Given $x\in \overline{\cup_{k\in\mathbb{Z}} L_k}$, let $(x_n)$ be a sequence in  $\cup_{k\in\mathbb{Z}} L_k$ with  $x_n\underset{n\to\infty}{\to }x$. Let $\epsilon$ with
\[
0<\epsilon <\mbox{min}\{ 1, \mbox{dist}(x, \partial\Omega )\}.
\]
Since $\psi$ is strongly run-away there exists $r\in\mathbb{N}$ so that for each $n\ge rN$ we have
\[
\psi_n(\overline{B(x,\epsilon)}\cup L)\cap (\overline{B(x,\epsilon)}\cup L) =\emptyset.
\]
Thus for each $n\ge r$ we must have
\begin{equation}\label{eq:iii2}
L_n\cap \overline{B(x, \epsilon)}=\psi_{nN}(L)\cap \overline{B(x, \epsilon)} \subseteq \psi_{nN}(\overline{B(x,\epsilon)}\cup L)\cap (\overline{B(x,\epsilon)}\cup L) =\emptyset
\end{equation}
and 
\begin{equation} \label{eq:iii3}
\overline{B(x, \epsilon)}\cap \psi_{nN}^{-1}(L)\subseteq  (\overline{B(x,\epsilon)}\cup L) \cap \psi_{nN}^{-1} (\overline{B(x,\epsilon)}\cup L)=\emptyset.
\end{equation}
By $\eqref{eq:iii3}$ for each $n\ge r$ we have
\begin{equation} \label{eq:iii4}
\overline{B(x,\epsilon)}\cap L_{-n}=\emptyset,
\end{equation}
so by $\eqref{eq:iii2}$ and $\eqref{eq:iii4}$ the tail of the sequence $(x_n)$ must lie in $\cup_{-r\le k\le r} L_k$, giving
\[
x\in \overline{\cup_{-r\le k\le r} L_k}=\cup_{-r\le k\le r} L_k\subseteq \cup_{k\in\mathbb{Z}}L_k.
\]
So $\cup_{k\in\mathbb{Z}}L_k$ is closed, and by $(ii)$ we have 
\[
\begin{aligned}
\partial (\cup_{k\in \mathbb{Z}}L_k)  &\subseteq \cup_{k\in\mathbb{Z}} L_k \setminus \mbox{int}(\cup_{k\in\mathbb{Z}} L_k) \\
&= \cup_{k\in\mathbb{Z}} L_k \setminus \cup_{k\in\mathbb{Z}}\mbox{int}( L_k) \\
&=\cup_{k\in\mathbb{Z}} (L_k\setminus \mbox{int}(L_k) )\\
&=\cup_{k\in\mathbb{Z}} \partial L_k.
\end{aligned}
\]
But the inclusion $\partial  (\cup_{k\in \mathbb{Z}}L_k)   \supseteq \cup_{k\in\mathbb{Z}} \partial L_k$ easily follows by $(i)$.
\end{proof}
The proof of Theorem~\ref{T:4.2} is now complete.
\end{proof}

\section{The one-dimensional case} 

\begin{theorem} \label{T:2.2} Let $\emptyset\ne \Omega \subset \mathbb{R}$ be open, and let $\psi \in C^{p} (\Omega, \Omega)$ and $\omega \in C^{p} (\Omega, \mathbb{K})$, where $p\in\mathbb{N}\cup\{\infty\}$. 
Then $C_{\omega, \psi}$ is supercyclic on $C^{p} (\Omega, \mathbb{K})$ if and only if  it is mixing, and if and only if it is chaotic.
\end{theorem}

To prove Theorem~\ref{T:1.3}, Przestacki showed that for a continuous and injective self-map $\psi$ of $\mathbb{R}$, having no fixed points is equivalent to being strongly run-away. 
To prove Theorem~\ref{T:2.2} we extend this result to any open subset $\Omega$ of $\mathbb{R}$.

\begin{lemma} \label{L:run-away,interval} Let $\psi: (a, b) \to (a, b)$ be continuous and injective where $-\infty \leq a < b \leq \infty$, then the following are equivalent:
\begin{enumerate}
\item[(1)]\ $\psi$ has no fixed points.
\item[(2)]\ $\psi$ is run-away.
\item[(3)]\ $\psi$ is strongly run-away.
\end{enumerate}
\end{lemma}

\begin{proof}
The implications $(3) \Rightarrow (2)$ and $(2) \Rightarrow(1)$ are immediate, so we show $(1) \Rightarrow (3)$. The case where $(a, b) = (-\infty, \infty)$  is \cite[Lemma 4.1]{Przestacki2017}. Let's start with the case where a and b are both finite. Since $\psi$ is fixed point free, we either have $\psi(x) > x$ or $\psi(x) < x$ for all $x \in (a,b)$. Without loss of generality, we assume $\psi(x) > x$  for all $x \in (a,b)$, with the other case being similar. Since $\psi$ is injective and has no fixed points, the Intermediate Value Property now ensures that $\psi$  is strictly increasing. So for each $a < c < d < b$ we have
\[
\psi([c,d]) = [\psi(c), \psi(d)].
\]
Let $K \subset (a, b)$ be compact.  It suffices to find an $N \in \mathbb{N}$ such that $\psi_n(K) \cap K = \emptyset$ for $n \geq N$. Enlarging K if necessary, we may assume that $K = [c,d]$ for some $\emptyset\ne [c, d]\subset (a, b)$.

Notice that the sequence $(\psi_n (c))$ is increasing, so $\lim_{n\to\infty} \psi_n (c)=b$. Otherwise there exists $a_0 \in (a,b)$ such that $(\psi_n (c))$ converges to $a_0$ and the continuity of  $\psi$ forces

\[
a_0 = \lim\limits_{n \to\ \infty} \psi_n (c)\ = \lim\limits_{n \to\ \infty} \psi_{n+1} (c)\ = \lim\limits_{n \to\ \infty} \psi(\psi_n (c)) = \psi( \lim\limits_{n \to\ \infty} \psi_n (c)) = \psi(a_0),\]
 contradicting that $\psi$ has no fixed points. So $\lim_{n\to\infty} \psi_n (c)=b$ and thus  there exists $N \in \mathbb{N}$ such that $\psi_n (c) > d$ for each $n \geq N$ and  thus
\[
\psi_n(K) \cap K = \psi_n([c, d]) \cap [c,d] = [\psi_n (c), \psi_n (d)] \cap [c, d] = \emptyset
\]
 for each $n \geq N$. So $\psi$ is strongly run-away. The proofs of the cases where $a = -\infty$ or $b = \infty$ have one difference: the monotone sequence $(\psi_n(c))$ may be unbounded. But the method to show that $\psi$ is strongly run-away remains the same. 
\end{proof}

We next consider with Lemma~\ref{L:run-away,finite} below the case when $\Omega$ is finitely connected. We first recall the following elementary fact.

\begin{lemma} \label{L:subseq,strongly} Let 
$\psi$  be a  continuous self-map on a topological space $\Omega$. The following are equivalent:
\begin{enumerate}
\item[(i)] $\psi $ is strongly run-away.
\item[(ii)]\ For each $p\in\mathbb{N}$, the $p$-th iterate $\psi_p$ of $\psi$   is strongly run-away.
\item[(iii)]\ For some $p\in\mathbb{N}$, the $p$-th iterate $\psi_p$ of $\psi$   is strongly run-away.
\end{enumerate}
\end{lemma}
\begin{proof} The implications $(i)\Rightarrow (ii)$ and $(ii)\Rightarrow (iii)$ are immediate. To see $(iii)\Rightarrow (i)$, fix a compact set $K \subset \Omega$ and consider the compact set $K_1 := K \cup \psi(K) \cup ... \cup \psi_{p-1} (K)$. Since $\psi_p$ is strongly run-away, there exists $N \in \mathbb{N}$ such that for each $n \geq N$ we have
\[
\psi_{pn} (K_1) \cap K_1 = \emptyset,\]
 so for each $n\ge N$ and $j=1,\dots , p-1$ we have
 \[
 K\cap \psi_{pn-j}(K)=\emptyset,
 \]
 and thus $K\cap \psi_m(K)=\emptyset $ for each $m\ge pN$.
\end{proof}

\begin{lemma} \label{L:run-away,finite} Let $\Omega = I_1 \cup ... \cup I_m$, where $I_1,\dots , I_m$ are pairwise disjoint non-empty open intervals, and let $\psi: \Omega \to \Omega$ be injective and continuous. Then the following are equivalent:
\begin{enumerate}
\item[(1)]\ $\psi$ has no fixed points, no points of period 2,\dots , and not points of period $m$.
\item[(2)]\ $\psi$ is run-away.
\item[(3)] $\psi$ is strongly run-away.
\end{enumerate}
\end{lemma}

\begin{proof} 
The implications $(3) \Rightarrow (2)$ and $(2) \Rightarrow (1)$ are immediate. We show  $(1)\Rightarrow (3)$ by induction on $m$. The case $m=1$ is Lemma~\ref{L:run-away,interval}. For the inductive step, assume $2\le m$. Notice  that by a connectedness argument for each $j=1,\dots ,m$ there exists a unique $\sigma(j)\in \{ 1,\dots, m\}$ so that
$
\psi(I_j)\subseteq I_{\sigma(j)}$,
so we have two cases:

\noindent Case 1: There exists $j\in \{1,\dots, m\}$ so that $\psi(\Omega)\subseteq \Omega\setminus I_j$. Here $\psi$ is strongly run-away on $\Omega\setminus I_j$ by the inductive assumption, and since $\psi(I_j)\subseteq \Omega\setminus I_j$ it follows that $\psi$ is strongly run-away on $\Omega$.

\noindent Case 2: $\psi (\Omega) \cap I_j \neq \emptyset$ for each $j=1,\dots, m$. In this case $\psi$ acts as a permutation on the indices of the m intervals. That is, there exists a bijection $\sigma:\{ 1,2,\dots, m\}\to \{ 1,2,\dots, m\}$ so that 
\[
\psi(I_j)\subseteq I_{\sigma(j)}
\]
for each $j=1,\dots,m$.  Hence if $p$ is the order of $\sigma$ in the group of permutations of $m$ elements,  then  $\psi_p(I_j)\subseteq I_j$ for each $j=1,\dots m$. So $\psi_p$ is strongly run-away on each $I_j$ by Lemma \ref{L:run-away,interval}, and thus  $\psi_p$ is strongly run-away on $\Omega$. So $\psi$ is strongly run-away on $\Omega$ by Lemma~\ref{L:subseq,strongly}. 
\end{proof}

Finally, we consider the infinitely connected case.

\begin{lemma} \label{L:run-away,infinite} Let $\Omega = \cup_{k=1}^{\infty} I_k$ where $\{I_k\}$ is a collection of pairwise disjoint non-empty open intervals, and let $\psi: \Omega \to \Omega$ be injective and continuous. The following are equivalent:
\begin{enumerate}
\item[(1)]\ $\psi$ has no periodic points.
\item[(2)]\ $\psi$ is run-away.
\item[(3)]\ $\psi$ is strongly run-away. 
\end{enumerate}
\end{lemma}

\begin{proof} Once again, we only have to show that (1) implies (3). The key is to rewrite $\Omega$ into a union of pairwise disjoint $\psi$-invariant sets that are minimal in size. For each $\emptyset\ne A\subset\mathbb{N}$ we define
\begin{equation} \label{eq:VA}
V_A:=\cup_{j\in A} I_j,
\end{equation}
and for each $k\in \mathbb{N}$ define
\[
\Omega_k:=\cap_{\{ A\subseteq\mathbb{N}: \ k\in A \mbox{ and } \psi(V_A)\subseteq V_A \} } V_A.
\] 
Notice that $I_k\subset \Omega_k$, and that $\psi(\Omega_k)\subseteq \Omega_k$. We refer to  $\Omega_k$ as the Invariant Component of $\Omega$ that contains $I_k$. We see next that each Invariant Component is of the form~$\eqref{eq:VA}$.

\noindent Claim 1: For each $k\in \mathbb{N}$ there exists a unique  $A_k^*\subseteq \mathbb{N}$ so that 
\[
\Omega_k=V_{A_k^*}.
\]
To see the claim, let $x\in \Omega_k$. So $x\in I_\ell$ for a unique $\ell\in\mathbb{N}$. It suffices to show that $I_\ell\subseteq \Omega_k$.
Let $A$ be a subset of $\mathbb{N}$ with $k\in A$ and $\psi(V_A)\subseteq V_A$. We want to show that $I_\ell \subseteq V_A$.
Since $x\in \Omega_k\subseteq V_A$, there exists $j\in A$ so that $x\in I_j$. Hence since $I_j$ and $I_\ell$ are connected components of $\Omega$, we must have 
\[
I_\ell=I_j\subseteq V_A
\]
and Claim 1 holds.  Notice that by Claim 1 for each $k,\ell\in \mathbb{N}$ we must have 
\[
\Omega_k\cap\Omega_\ell=\emptyset \mbox{ or  } \Omega_k=\Omega_\ell.\]

We next show that $\psi$ is strongly run-away on each invariant component:

\noindent Claim 2: Let $k\in \mathbb{N}$. Then $\psi$ is strongly run-away on $\Omega_k$.

To see Claim 2, let $A_k^*\subseteq \mathbb{N}$ so that 
\[
\Omega_k=V_{A_k^*}=\cup_{j\in A_k^*} I_j.
\]  
If $A_k^*$ is finite then $\Omega_k$ is finitely connected and Claim 2 follows from Lemma~\ref{L:run-away,finite}. So assume that $A_k^*$ is infinite. Notice that in this case for each $j\in A_k^*$ and $n\in\mathbb{N}$ each of the intervals
\begin{equation} \label{eq:FI2}
\psi(I_j), \psi_2(I_j),\dots, \psi_n(I_j)
\end{equation}
must be contained in a different connected component of $\Omega_k=V_{A_k^*}=\cup_{j\in A_k^*} I_j$.  Now, if $K\subset \Omega_k$ is compact there exist $j_1,\dots, j_r\in A_k^*$ so that
\[
K\subset \cup_{\ell=1}^r I_{j_\ell}.
\]
By $\eqref{eq:FI2}$, for each $\ell=1,\dots, r$ there exists $N_\ell\in\mathbb{N}$ so that
\[
\psi_n(I_{j_\ell})\cap \cup_{\ell=1}^r I_{j_\ell}=\emptyset
\]
for each $n\ge N_\ell$. Hence 
\[
\psi_n(K)\cap K=\emptyset 
\]
for each $n\ge \max{\{ N_1,\dots N_r\} }$, and Claim 2 holds.

To see that $\psi$ is strongly run-away on $\Omega$, let $K\subseteq \Omega$ be compact. Since the collection of invariant components of $\Omega$ is an open cover of $K$, there exists $r\in \mathbb{N}$ so that 
\[
K\subset \cup_{\ell=1}^r \Omega_\ell.
\]
Since $\psi$ is strongly run-away on each $\Omega_\ell$ and $\psi(\Omega_\ell)\subseteq \Omega_\ell$, there exists $N\in\mathbb{N}$ so that 
for each $1\le \ell\le r$ and $n\ge N$ we have
\[
\psi_n(K\cap \Omega_\ell)\cap (K\cap \Omega_\ell) =\emptyset,
\]
which forces
\[
\psi_n(K)\cap K=\emptyset
\]
for each $n\ge N$.
\end{proof}

By Lemma~\ref{L:run-away,interval}, Lemma~\ref{L:run-away,finite} and Lemma~\ref{L:run-away,infinite} we have the following.

\begin{corollary} \label{C:runaway=strongrunaway} Let $\Omega \subseteq \mathbb{R}$ be open. If $\psi: \Omega \to \Omega$ is injective and continuous, then it is run-away if and only if it is strongly run-away and if and only if it has no periodic points.
\end{corollary}
\begin{remark}
While lacking periodic points ensures  an injective continuous map $\psi:\Omega\to \Omega$ to be strongly run-away in the one-dimensional setting, this is no longer true in higher dimensions. Indeed, consider the open set $\Omega = \mathbb{R}^{2} \setminus \{ (0,0) \}$. Let $\psi$ be a rotation by an angle $\theta$ that is an irrational multiple of $\pi$. So $\psi$ is a self map of $\Omega$ that is injective and continuous. Moreover, given $x \in \Omega$, by the selection of $\theta$ the orbit of x under $\psi$ must be dense in the circle of radius $|x|$ that is centered at the origin. That is, $\psi$ has no periodic points and $\psi$ is not run-away since circles centered at the origin are compact and invariant under $\psi$.
\end{remark}

We are ready to show Theorem \ref{T:2.2}.

\begin{proof}[Proof of Theorem \ref{T:2.2}] Each mixing operator is  supercyclic, so this happens for $C_{\omega, \psi}$. Conversely, suppose that $C_{\omega, \psi}$ is supercyclic. 
By Theorem \ref{T:2.1}, $\psi$ is run-away and Conditions (3) (i)-(iii) of Theorem~\ref{T:2.1} hold. Since $\Omega\subseteq \mathbb{R}$, Corollary~\ref{C:runaway=strongrunaway} ensures that 
$\psi$ is strongly run-away. So  $C_{\omega, \psi}$ is mixing by Theorem~\ref{T:2.1}. For the equivalence with chaos, recall  that every chaotic operator is supercyclic and the converse of this now follows from Theorem~\ref{T:4.2}, as we already established that $C_{\omega, \psi}$ is mixing whenever it is supercyclic and $\Omega\subseteq \mathbb{R}$.
\end{proof}

\section{Final Comments}
Since in general for any operator on a separable Fr\'echet space we have the implications
\[
\mbox{mixing}\implies \mbox{weak-mixing}\implies \mbox{transitive}\implies \mbox{supercyclic}
\]
and
\[
\mbox{chaos}\implies \mbox{weak mixing},
\]
it is natural to seek when any of the converses of these implications hold within a certain class of operators. For a weighted composition operator $C_{\omega, \psi}$ acting on the space $H(\Omega )$ of holomorphic functions on an arbitrary domain $\Omega$ of the complex plane, for instance, Golinski and Przestacki~\cite{golprz2021} fully characterized hypercyclicity and showed that $C_{\omega, \psi}$ is hypercyclic if and only if it is weakly mixing. Indeed, the latter is known to hold if and only if $C_{\omega, \psi}$ is mixing, and except possibly when $\Omega$ is conformally equivalent to a punctured disc the operator  $C_{\omega, \psi}$ is supercyclic on $H(\Omega)$ if and only if it is mixing \cite{besfoster2024}.

Motivated by the above equivalences holding for weighted composition operators acting on $H(\Omega)$ and the corresponding given equivalences we saw with Theorem~\ref{T:2.2} on $C^p(\Omega, \mathbb{K})$ when $\Omega\subseteq \mathbb{R}$ it is natural to ask the following.

\begin{question} \label{P1}
Let  $d \geq 2$ and $\Omega \subset \mathbb{R}^{d}$ be open. Is every supercyclic weighted composition operator on $C^{\infty} (\Omega, \mathbb{K})$ also mixing? Equivalently, is every injective smooth function $\psi:\Omega \to \Omega$ strongly run-away whenever it is run-away?
\end{question}

With Theorem~\ref{T:2.2} and Theorem~\ref{T:2.1} we have a characterization for chaos of weighted composition operators on $C^{\infty} (\Omega, \mathbb{K})$  whenever $\Omega$ is an open subset of $\mathbb{R}$. Hence Problem~\ref{P2} may be reformulated as follows.
 
 \setcounter{problem}{1} 
\begin{problem} \label{Q3}
Let $\Omega\subseteq \mathbb{R}^d$ be open. Is there a characterization for chaos for weighted composition operators on  $C^{\infty} (\Omega, \mathbb{K})$ when $d\ge 2$?
\end{problem}

Theorem~\ref{T:2.2} also motivates the following.

\begin{question} \label{P2}
Let  $d \geq 2$ and $\Omega \subset \mathbb{R}^{d}$ be open. Is every supercyclic weighted composition operator on $C^{\infty} (\Omega, \mathbb{K})$ also chaotic?
\end{question}
\noindent
By Theorem~\ref{T:4.2}, a positive answer to Question~\ref{P1} gives a positive answer to Question~\ref{P2}.
Finally, motivated again by Theorem~\ref{T:2.2}  we also ask whether the converse of Theorem~\ref{T:4.2} holds in the multidimensional case.

\begin{question} \label{P3}Let  $d \geq 2$ and $\Omega \subset \mathbb{R}^{d}$ be open.
Is every chaotic weighted composition operator on $C^{\infty} (\Omega, \mathbb{K})$ also mixing?
\end{question}

\noindent
By Theorem~\ref{T:2.1} and Theorem~\ref{T:4.2},
a positive answer to either  Question~\ref{P2} or Question~\ref{P3} would solve Problem~\ref{Q3}.

\end{document}